\documentclass[11pt]{article}
\usepackage{amsfonts}
\usepackage{graphics}
\usepackage{color}
\usepackage{cite}
\usepackage{latexsym}
\usepackage[paper=a4paper, left=2.4cm, right=2.4cm, top=2.5cm, bottom=1.8cm, headheight=5.5pt, footskip=1.0cm, footnotesep=0.8cm, centering, includefoot]{geometry}
\usepackage{amsmath}
\allowdisplaybreaks
\usepackage{amssymb}
\usepackage[dvips]{epsfig}
\usepackage{amscd}

\newtheorem{theorem}{Theorem}[section]
\newtheorem{remark}{Remark}[section]
\newtheorem{definition}{Definition}[section]
\newtheorem{lemma}[theorem]{Lemma}

\newcommand{\vp}{\varphi}
\newcommand{\n}{\rho}

\newcommand{\ti}{\tilde}

\DeclareMathOperator{\loc}{loc}

\renewcommand{\div}{ {\rm div }  }

\newcommand{\pa}{\partial}
\renewcommand{\r}{\mathbb{R}}

\newcommand{\bt}{\begin{theorem}}
\newcommand{\bl}{\begin{lemma}}
\newcommand{\el}{\end{lemma}}
\newcommand{\et}{\end{theorem}}

\newcommand{\al}{\alpha}

\newcommand{\ve}{\varepsilon}
\newcommand{\la}{\label}

\newcommand{\bn}{\begin{eqnarray}}
\newcommand{\en}{\end{eqnarray}}
\newcommand{\bnn}{\begin{eqnarray*}}
\newcommand{\enn}{\end{eqnarray*}}

\newcommand{\bnnn}{\begin{eqnarray*}}
\newcommand{\ennn}{\end{eqnarray*}}
\newcommand{\ba}{\begin{aligned}}
\newcommand{\ea}{\end{aligned}}
\newcommand{\be}{\begin{equation}}
\newcommand{\ee}{\end{equation}}
\def\O{{\r^2 }}
\def\p{\partial}
\def\norm[#1]#2{\|#2\|_{#1}}

\def\la{\label}

\def\na{\nabla}

\makeatletter
\@addtoreset{equation}{section}
\makeatother

\makeatletter
\@addtoreset{equation}{section}
\makeatother

\title{Global Existence and Large Time Asymptotic Behavior of Strong Solutions to the Cauchy Problem of 2D Density-Dependent Navier-Stokes Equations with Vacuum\thanks{B. L\"u is supported by NNSFC Tianyuan (No. 11426131) and Natural Science Foundation of Jiangxi Province (No. 20142BAB211006). X. Shi is supported by  NNSFC (Nos.  11371348 \& 11471321).}
}
\author{ Boqiang L\"u\thanks{College of Mathematics and Information Science, Nanchang Hangkong University, Nanchang 330063,
People's Republic of China ({\tt lvbq86@163.com}). }
\quad Xiaoding Shi\thanks{Department of Mathematics, School of Science, Beijing University of Chemical Technology, Beijing 100029,
People's Republic of China ({\tt shixd@mail.buct.edu.cn}).}
\quad Xin Zhong\thanks{Corresponding author. Institute of Applied Mathematics, AMSS, Chinese Academy of Sciences, Beijing 100190,
People's Republic of China ({\tt xzhong1014@amss.ac.cn}).
 }
 }

\date{ }

\begin{document}
\maketitle

\begin{abstract}
We are concerned with the Cauchy problem of the two-dimensional (2D) nonhomogeneous incompressible Navier-Stokes equations with vacuum as far-field density. It is proved that if the initial density decays not too slow at infinity, the 2D Cauchy problem of the density-dependent Navier-Stokes equations on the whole space $\mathbb{R}^2$ admits a unique global strong solution. Note that the initial data can be arbitrarily large and the initial density can contain vacuum states and even have compact support. Furthermore, we also obtain the large time decay rates of the spatial gradients of the velocity and the pressure which are the same as those of the homogeneous case.
\end{abstract}


Math Subject Classification: 35Q35; 76N10.

\section{Introduction}
The motion of a two-dimensional nonhomogeneous incompressible fluid is governed by the following Navier-Stokes equations:
\be\label{1.1}
\begin{cases}
\partial_{t}\rho+\div(\rho u)=0, \\
\partial_{t}(\rho u )+\div(\rho u \otimes u )-\mu\Delta u +\nabla P=0, \\
\div u=0.
\end{cases}
\ee
Here, $t\ge 0$ is time, $x\in\r^2$ is the spatial coordinates, and the unknown functions $\rho=\rho(x,t)$, $u=(u^1,u^2)(x,t)$, and $P=P(x,t)$ denote the density, velocity, and pressure of the fluid, respectively;
$\mu>0$ stands for the viscosity constant.

We consider the Cauchy problem of \eqref{1.1} with $(\rho, u )$ vanishing at infinity (in some weak sense) and the initial conditions:
\begin{equation}\label{1.3}
\rho(x,0)=\rho_0(x),\ \rho u (x,0)=\rho_0 u _0(x), \ x\in\mathbb{R}^2,
\end{equation}
for given initial data $\rho_0$ and $u_0$.

The study of the system \eqref{1.1} has a long history. Since the pioneering work of Leray \cite{L1934}, there are huge literatures on the studies of the large time existence and behavior of solutions of \eqref{1.1} in the case that the density is constant (refer to \cite{CF1988,K1984,L1969,L1996,T2001} and references therein). The mathematical study of nonhomogeneous incompressible flow goes back to the seventies of the last century.
When the initial density is strictly away from vacuum, Kazhikov \cite{K1974} established the global existence of weak solutions (see also \cite{AK1973}). Later, Antontsev-Kazhikov-Monakhov \cite{AKM1990} gave the first result on local existence and uniqueness of strong solutions. Moreover, they also proved that the unique local strong solution is a global one in two dimensions.
Recently, the global existence result of \eqref{1.1} with small initial data belonging to certain scale invariant spaces was obtained in \cite{GZ2009}.

On the other hand, when the initial density allows vacuum, Simon \cite{S1990} proved the global existence of weak solutions, which was extended by Lions \cite{L1996} to the case of density-dependent viscosity. For the initial density allowing vacuum, Choe-Kim \cite{CK2003} proposed a compatibility condition and obtained the local existence of strong solutions of \eqref{1.1} with large data in two-dimensional bounded domains. Recently, Huang-Wang \cite{HW2014} showed that the local strong solution obtained in \cite{CK2003} is indeed a global one. Very recently, motivated by \cite{LL2014,lh1}, L\"u-Xu-Zhong \cite{LZ2015} establihsed the local (with generally large data) existence of strong solutions to 2D Cauchy problem of the nonhomogeneous incompressible magnetohydrodynamic equations with vacuum as far-field density, which in particular yields the local existence of strong solutions to the 2D Cauchy problem of \eqref{1.1}. However, the global existence of strong solution to the 2D Cauchy problem of \eqref{1.1} with vacuum and general initial data is still open. In fact, this is the main aim in this paper.

Before stating the main results, we first explain the notations and
conventions used throughout this paper. For $R>0$, set
\begin{equation*}
B_R \triangleq\left.\left\{x\in\r^2\right|\,|x|<R \right\},
\quad \int \cdot dx\triangleq\int_{\r^2}\cdot dx.
\end{equation*}
Moreover, for $1\le r\le \infty$ and $k\ge 1$, the standard Lebesgue and  Sobolev spaces are defined as follows:
\bnn
L^r=L^r(\r^2 ),\quad W^{k,r}= W^{k,r}(\r^2), \quad H^k = W^{k,2}.
\enn

Now we define precisely what we mean by strong solutions.
\begin{definition}\label{def1}
If all derivatives involved in \eqref{1.1} for $(\rho,u,P)$ are regular distributions, and equations \eqref{1.1} hold almost everywhere in $\mathbb{R}^2\times(0,T)$, then $(\rho,u,P)$ is called a strong solution to \eqref{1.1}.
\end{definition}

Without loss of generality, we assume that the initial density $\n_0$ satisfies
\be\la{oy3.7}
\int_{\r^2} \n_0dx=1,
\ee
which implies that there exists a positive constant $N_0$ such that
\be\label{13}
\int_{B_{N_0}}  \n_0  dx\ge \frac12\int\n_0dx=\frac12.
\ee

Our main result can be stated as follows:
\begin{theorem}\label{thm1}
In addition to \eqref{oy3.7} and \eqref{13}, assume that the initial data $(\rho_0, u _0)$ satisfy for any given numbers $  a>1$ and $q>2$,
\begin{equation}\label{2.2}
\rho_{0}\ge 0,\  \rho_{0}\bar{x}^{a}\in L^{1}\cap H^{1}\cap W^{1,q},\ \div u _0=0,\ \nabla u _{0}\in L^2,\  \sqrt{\rho_0} u _0\in L^2,
\end{equation}
where
\begin{equation}\label{2.1}
\bar{x}\triangleq(e+|x|^2)^{1/2}\log^{2}(e+|x|^2).
\end{equation}
Then  the problem \eqref{1.1}-\eqref{1.3} has a unique global strong solution $(\rho, u, P)$ satisfying that
for any $0<T<\infty$,
\begin{equation}\label{2.3}
\begin{cases}
0\le \rho\in C([0,T];L^{1}\cap H^{1}\cap W^{1,q}), &\\
\rho\bar{x}^{a}\in L^{\infty}(0,T;L^{1}\cap H^{1}\cap W^{1,q}), &\\
\sqrt{\rho} u ,\ \nabla u ,\ \bar{x}^{-1} u ,\
\sqrt{t}\sqrt{\rho} u_t,~ \sqrt{t}\na P, ~\sqrt{t}\na^2 u \in L^{\infty}(0,T;L^2), &\\
\nabla u \in L^{2}(0,T;H^1)\cap L^{(q+1)/q}(0,T;W^{1,q}), &\\
\na P\in  L^2(0,T;L^2)\cap  L^{(q+1)/q}(0,T;L^q), \\
\sqrt{t}\nabla u \in L^{2}(0,T;W^{1,q}), &\\
\sqrt{\rho}u_t ,\ \sqrt{t}\nabla u_t ,\ \sqrt{t}\bar{x}^{-1}u_t\in L^{2}(\mathbb{R}^{2}\times(0,T)),
\end{cases}
\end{equation}
and
\begin{equation}\label{d}
\inf_{0\leq t\leq T}\int_{B_{N_{1}}}\rho(x,t)dx\geq\frac{1}{4},
\end{equation}
for some positive constant $N_1$ depending only on $\|\sqrt{\n_0} u_0\|_{L^2}, N_0$, and $T$. Moreover, $(\rho,u,P)$ has the following decay rates, that is, for $t\geq1$,
\be \la{e}
\|\na u(\cdot,t)\|^2_{L^2}+
\|\na^2 u(\cdot,t)\|_{L^2}+\|\na P(\cdot,t)\|_{L^2}\le Ct^{-1},
\ee
where $C$ depends only on $\mu, \|\n_0\|_{L^1\cap L^\infty}, \|\sqrt{\n_0} u_0\|_{L^2}$,  and $\|\na u_0\|_{L^2} $.
\end{theorem}

\begin{remark}
Our Theorem \ref{thm1} holds for arbitrarily large initial data which is in sharp contrast to Li-Xin \cite{LX2014} where the smallness conditions on the initial energy is needed in order to obtain the global existence of strong solutions to the 2D compressible Navier-Stokes equations.
\end{remark}

\begin{remark}
It should be noted here that although the equations \eqref{1.1}$_2$ degenerate near vacuum, our large time decay rates \eqref{e} are the same as those of the homogeneous case \cite{K1984}.
\end{remark}

\begin{remark}
Compared with \cite{HW2014}, there is no need to impose the additional compatibility conditions on the initial data for the global existence of the strong solution.
\end{remark}

We now make some comments on the key ingredients of the analysis in this paper. Note that for initial data satisfying \eqref{2.2}, L\"u-Xu-Zhong \cite{LZ2015} established the local existence and uniqueness of strong solutions to the Cauchy problem \eqref{1.1}-\eqref{1.3} (see Lemma \ref{lem21}). Thus, to extend the local strong solution to be a global one, one needs  to obtain global a priori estimates on strong solutions to \eqref{1.1}-\eqref{1.3} in suitable higher norms. It should be pointed out that the crucial techniques of proofs in \cite{HW2014} cannot be adapted directly to the situation treated here, since their arguments depend crucially on the boundedness of the domains. Moreover, it seems difficult to bound the $L^q(\mathbb{R}^2)$-norm of $ u $ just in terms of $\|\sqrt{\rho} u \|_{L^{2}(\mathbb{R}^2)}$ and $\|\nabla u \|_{L^{2}(\mathbb{R}^2)}$.
To overcome these difficulties mentioned above, some new ideas are needed.
First, motivated by \cite{H19951,HLX2012,LX2014}, multiplying \eqref{1.1}$_2$ by the material derivatives of the velocity, $\dot u\triangleq u_t+u\cdot\na u$, instead of the usual $u_t$ (see \cite{HW2014}), we find that the key point  to obtain the estimate on the $L^\infty(0,T;L^2(\r^2))$-norm of the gradient of the velocity is to bound the term
\begin{equation*}
I_2\triangleq\int P\pa_iu^j\pa_ju^i dx.
\end{equation*}
Motivated by \cite{D1997}, we find that $I_2$ in fact can be bounded by  $ \|\na P\|_{L^2}\|\na u\|_{L^2}^2$ (see \eqref{3.8}) since $\pa_iu^j\pa_ju^i\in \mathcal{H}^{1}$ due to the facts that $\div u=0 $ and that $\na^\perp \cdot \na u=0$ (see Lemma \ref{lem27} also). Next, to obtain the estimates on the gradient of the density, motivated by \cite{H19951,HL20130,HL20131,LX2014},   we apply the operator $\pa_t+ u\cdot \na $ to \eqref{1.1}$_2^j$ and multiply the resultant equality by $\dot u^j$ to get the time-independent estimates on both the $L^\infty(0,T;L^2(\r^2))$-norm of $t^{1/2} \n^{1/2}\dot u$ and the $L^2(\r^2\times(0,T ))$-norm of $t^{1/2}\na \dot u$ (see \eqref{133}) which combined with some careful analysis on the spatial weighted estimate of the density (see \eqref{06.1}) thus yield the bound on the $L^1(0,T;L^\infty)$-norm of the gradient of the velocity (see \eqref{nu1}) which in particular implies the desired bound on the $L^\infty(0,T;L^q)$-norm of the gradient of the density. Finally, with these a priori estimates on
the gradients of the density and the velocity at hand, one can
estimate the higher order derivatives by using the same arguments
as in \cite{LX2014,HW2014} to obtain the desired results.

The rest of this paper is organized as follows. In Section \ref{sec2}, we collect some elementary facts and inequalities that will be used later. Section \ref{sec3} is devoted to the a priori estimates. Finally, we will give the proof of Theorem \ref{thm1} in Section \ref{sec4}.

\section{Preliminaries}\label{sec2}
In this section we shall enumerate some auxiliary lemmas used in this paper.

We start with the local existence of strong solutions whose proof can be found in \cite[Theorem 1.2]{LZ2015}.
\begin{lemma}\label{lem21}
Assume that $(\rho_0, u _0)$ satisfies \eqref{2.2}. Then there exist a small time $T>0$ and a unique strong solution $(\rho, u, P)$ to the problem \eqref{1.1}-\eqref{1.3} in $\mathbb{R}^{2}\times(0,T)$ satisfying \eqref{2.3} and \eqref{d}.
\end{lemma}

Next, the following Gagliardo-Nirenberg inequality (see \cite{N1959}) will be used later.
\begin{lemma}[Gagliardo-Nirenberg]\label{lem22}
For $q\in[2,\infty), r\in(2,\infty)$, and $s\in(1,\infty)$, there exists some generic constant $C>0$ which may depend on $q,$ $r$, and $s$ such that for $f\in H^{1}(\mathbb{R}^2)$ and $g\in L^{s}(\mathbb{R}^2)\cap D^{1,r}(\mathbb{R}^2)$, we have
\begin{eqnarray*}
& & \|f\|_{L^q(\mathbb{R}^2)}^{q}\leq C\|f\|_{L^2(\mathbb{R}^2)}^{2}\|\nabla f\|_{L^2(\mathbb{R}^2)}^{q-2}, \\
& & \|g\|_{C(\overline{\mathbb{R}^2})}\leq C\|g\|_{L^s(\mathbb{R}^2)}^{s(r-2)/(2r+s(r-2))}\|\nabla g\|_{L^r(\mathbb{R}^2)}^{2r/(2r+s(r-2))}.
\end{eqnarray*}
\end{lemma}

The following weighted $L^m$ bounds for elements of the Hilbert space $\tilde{D}^{1,2}(\O)\triangleq\{ v \in H_{\loc}^{1}(\mathbb{R}^2)|\nabla v \in L^{2}(\mathbb{R}^2)\}$ can be found in \cite[Theorem B.1]{L1996}.
\begin{lemma} \la{1leo}
For $m\in [2,\infty)$ and $\theta\in (1+m/2,\infty),$ there exists a positive constant $C$ such that for all $v\in  \tilde{D}^{1,2}(\O),$ \be\la{3h} \left(\int_{\O} \frac{|v|^m}{e+|x|^2}l\left(\log \left(e+|x|^2\right)\right)^{-\theta}dx  \right)^{1/m}\le C\|v\|_{L^2(B_1)}+C\|\na v\|_{L^2(\O) }.\ee
\end{lemma}

The combination of Lemma \ref{1leo} and the Poincar\'e inequality yields
the following useful results on weighted bounds, whose proof can be found in \cite[Lemma 2.4]{LX2014}.

\begin{lemma}\label{lem26}
Let $\bar x$ be as in \eqref{2.1}. Assume that $\n \in L^1(\O)\cap L^\infty(\O)$ is a non-negative function such that
\bnn \la{2.i2}
\|\n\|_{L^1(B_{N_1})} \ge M_1, \quad \|\n\|_{L^1(\O)\cap L^\infty(\O)}\le M_2,
\enn
for positive constants $M_1, M_2$, and $ N_1\ge 1$. Then for $\ve> 0$ and $\eta>0,$ there is a positive constant $C$ depending only on $\ve,\eta, M_1,M_2$, and $ N_1$, such that every $v\in \ti D^{1,2}(\O)$ satisfies
\be\la{22}\ba
\|v\bar x^{-\eta}\|_{L^{(2+\ve)/\ti\eta}(\O)} &\le C \|\n^{1/2}v\|_{L^2(\O)}+C \|\na v\|_{L^2(\O)},
\ea\ee
with $\ti\eta=\min\{1,\eta\}$.
\end{lemma}

Finally, let $\mathcal{H}^{1}(\mathbb{R}^2)$ and BMO$(\mathbb{R}^2)$ stand for the usual Hardy and BMO spaces (see \cite[Chapter IV]{S1993}). Then the following well-known facts play a key role in the proof of Lemma \ref{lem3.3} in the next section.
\begin{lemma}\label{lem27}
(a) There is a positive constant $C$ such that
\begin{equation*}
\|E\cdot B\|_{\mathcal{H}^{1}(\O)}
\leq C\|E\|_{L^{2}(\O)}\|B\|_{L^{2}(\O)},
\end{equation*} for all $E\in L^{2}(\mathbb{R}^2)$ and $B\in L^{2}(\mathbb{R}^2)$ satisfying
\begin{equation*}
\div E=0,\ \nabla^{\bot}\cdot B=0\ \ \text{in}\ \ \mathcal{D}'(\mathbb{R}^2).
\end{equation*}
(b) There is a positive constant $C$ such that
\begin{equation}\label{lem1}
\| v \|_{{\rm BMO}(\O)}\leq C\|\nabla v \|_{L^{2}(\O)},
\end{equation} for all $  v \in \ti D^{1,2}(\mathbb{R}^2)$.
\end{lemma}
{\it Proof.}
(a) For the detailed proof, see \cite[Theorem II.1]{CLMS1993}.

(b) It follows  from the Poincar{\'e} inequality that for any ball $B\subset\mathbb{R}^2$
\begin{equation*}
\frac{1}{|B|}\int_{B}\left| v(x) - \frac{1}{|B|}\int_Bv(y)dy\right|dx\leq C\left(\int_{B}|\nabla v |^2dx\right)^{1/2},
\end{equation*}
which directly gives \eqref{lem1}.  \hfill $\Box$

\section{A Priori Estimates}\label{sec3}
In this section, we will establish some necessary a priori bounds for strong solutions $(\rho,u,P)$ to the Cauchy problem \eqref{1.1}-\eqref{1.3} to extend the local strong solution. Thus, let $T>0$ be a fixed time and $(\rho, u,P)$  be the strong solution to \eqref{1.1}-\eqref{1.3} on $\mathbb{R}^{2}\times(0,T]$ with initial data $(\rho_0,u_0)$ satisfying \eqref{oy3.7}-\eqref{2.2}.

In what follows, we will use the convention that $C$ denotes a generic positive constant depending on $\mu$, $a$, and the initial data, and use $C(\al)$ to emphasize that $C$ depends on $\al.$

\subsection{Lower Order Estimates}
First, since $\div u=0,$ we state the following well-known estimate on the $L^\infty(0,T;L^p)$-norm of the density.
\begin{lemma}[\cite{L1996}]\label{lem3.1}
There exists a positive constant $C$ depending only on $\|\n_0\|_{L^1\cap L^\infty}$ such that
\begin{equation}\label{3.1}
\sup_{t\in[0,T]}\|\rho\|_{L^{1}\cap L^\infty}\leq C.
\end{equation}
\end{lemma}

Next, the following lemma concerns the key time-independent estimates on the $L^\infty(0,T;L^2)$-norm of the gradient of the velocity.
\begin{lemma}\label{lem3.3}
There exists a positive constant $C$ depending only  on $\mu,   \|\n_0\|_{L^\infty}, \|\sqrt{\n_0} u_0\|_{L^2}$, and $\|\na u_0\|_{L^2}$
such that
\begin{equation}\label{3.5}
\sup_{t\in[0,T]} \|\nabla u \|_{L^2}^2 +\int_{0}^{T}\int\rho|\dot{ u
}|^{2}dxdt\leq C,
\end{equation}
where $\dot u\triangleq u_t+u\cdot\na u $ is the material derivatives of the velocity. Furthermore, one has
\begin{equation}\label{022}
\sup_{t\in[0,T]}t\|\nabla u \|_{L^2}^2
+\int_{0}^{T}t\int\rho|\dot{ u }|^{2}dxdt\leq C.
\end{equation}
\end{lemma}

{\it Proof.} First, applying standard energy estimate to \eqref{1.1} gives
\begin{equation}\label{0.11}
\sup_{t\in[0,T]}\|\sqrt{\rho} u \|_{L^2}^2
+\int_{0}^{T}\|\nabla u \|_{L^2}^{2}dt\leq C.
\end{equation}

Next, multiplying \eqref{1.1}$_2$ by $\dot{ u }$ and then integrating the resulting equality over $\mathbb{R}^2$ lead to
\begin{equation}\label{3.6}
\int\rho|\dot{ u }|^{2}dx=\int(\mu\Delta u \cdot\dot{ u }-\nabla P\cdot\dot{ u })dx\triangleq I_{1}+I_{2}.
\end{equation}
It follows from integration by parts and Garliardo-Nirenberg inequality that
\begin{align} \label{3.7}
I_{1}  &  =\int\mu\Delta u \cdot(\partial_{t} u + u \cdot\nabla u )dx \notag \\
 &  = -\frac{\mu}{2}\frac{d}{dt}\|\nabla u \|_{L^2}^{2}
-\mu\int\partial_{i}u^{j}\partial_{i}(u^{k}\partial_{k}u^{j})dx \notag \\
 &  \leq-\frac{\mu}{2}\frac{d}{dt}\|\nabla u \|_{L^2}^{2}+C\|\nabla u \|_{L^3}^{3} \notag \\
 &  \leq-\frac{\mu}{2}\frac{d}{dt}\|\nabla u \|_{L^2}^{2}
+C\|\nabla u \|_{L^2}^{2}\|\nabla^{2} u \|_{L^2}.
\end{align}
We deduce from integration by parts and \eqref{1.1}$_3$ that
\begin{align}\la{lz}
I_{2}  &  =-\int\nabla P\cdot(\partial_{t} u + u \cdot\nabla u )dx \notag \\
 &  = \int P\partial_{j}u^{i}\partial_{i}u^{j}dx \notag \\
  & \leq C\|P\|_{\rm BMO}\|\partial_{j}u^{i}\partial_{i}u^{j}\|_{\mathcal{H}^1},
\end{align}
where one has used the duality of $\mathcal{H}^1$ space and  BMO  one (see \cite[Chapter IV]{S1993}) in the last inequality. Since $\div(\partial_j u )=\partial_j\div u =0$ and $\nabla^{\bot}\cdot(\nabla u^{j})=0$,   Lemma \ref{lem27} yields
\begin{equation}\label{3.8}\ba
|I_{2}|=\left|\int P\partial_{j}u^{i}\partial_{i}u^{j}dx \right|\leq C\|\nabla P\|_{L^2}\|\nabla u \|_{L^2}^{2}.\ea
\end{equation}

Then, substituting \eqref{3.7} and \eqref{3.8} into \eqref{3.6} gives
\begin{equation}\label{3.9}
 \frac{\mu}{2}\frac{d}{dt}\|\nabla u \|_{L^2}^{2}+\|\sqrt{\rho}\dot{ u }\|_{L^2}^{2}
\leq C\left(\|\nabla^2 u \|_{L^2}+\|\nabla P\|_{L^2}\right)\|\nabla u \|_{L^2}^{2}.
\end{equation}
 On the other hand, since $(\n,u,P)$ satisfies the following Stokes system
\be\la{stokes1}
\begin{cases}
 -\mu\Delta u + \nabla P = -\n \dot u,\,\,\,\,&x\in \mathbb{R}^2,\\
 \div u=0,   \,\,\,&x\in  \mathbb{R}^2,\\
u(x)\rightarrow0,\,\,\,\,&|x|\rightarrow\infty,
\end{cases}
\ee
applying the standard $L^p$-estimate to \eqref{stokes1} (see \cite{T2001}) yields that for any $r\in (1,\infty),$
\be\label{001}\ba
\|\na^2 u \|_{L^r}+\|\nabla P\|_{L^r}& \le C(r)\|\n \dot u \|_{L^r}.
\ea\ee
It thus follows from \eqref{3.9}, \eqref{001}, and \eqref{3.1} that
\begin{eqnarray*}
\frac{\mu}{2}\frac{d}{dt}\|\nabla u \|_{L^2}^{2}+\|\sqrt{\rho}\dot{ u
}\|_{L^2}^{2} \leq C\|\sqrt{\rho}\dot{ u }\|_{L^2}\|\nabla u
\|_{L^2}^{2} \leq \frac{1}{2}\|\sqrt{\rho}\dot{ u
}\|_{L^2}^{2}+C\|\nabla u \|_{L^2}^{4} ,
\end{eqnarray*}
which implies that
\be\ba \label{3.11}
\frac{d}{dt}\left(\mu\|\nabla u \|_{L^2}^{2}\right)+\|\sqrt{\rho}\dot{ u
}\|_{L^2}^{2} \leq C\|\nabla u \|_{L^2}^{4} .
\ea\ee
This combined with \eqref{0.11} and Gronwall's inequality gives \eqref{3.5}.

Finally, multiplying \eqref{3.11} by $t$ leads to
\bnn\ba \label{3.12}
\frac{d}{dt}\left(t\mu\|\nabla u \|_{L^2}^{2}\right)+t\|\sqrt{\rho}\dot{ u }\|_{L^2}^{2}
\leq Ct\|\nabla u \|_{L^2}^{4} +C\|\nabla u \|_{L^2}^{2},
\ea\enn
which together with \eqref{0.11} and Gronwall's inequality yields \eqref{022} and completes the proof of Lemma \ref{lem3.3}.  \hfill $\Box$

Next, motivated by \cite{H19951,HLX2012,LX2014} where they deal with the compressible Navier-Stokes equations, we have the following estimates on the material derivatives of the velocity which are important for the higher order estimates of both the density and the velocity.

\begin{lemma}\label{lem3.4}
There exists a positive constant $C$ depending only on $\mu,   \|\n_0\|_{L^1\cap L^\infty}, \|\sqrt{\n_0} u_0\|_{L^2}$, and $\|\na u_0\|_{L^2} $ such that for $i=1,2$,
\begin{equation}\label{133}
\sup_{t\in[0,T]}t^i\|\sqrt{\rho}\dot{ u }\|_{L^2}^2
+\int_{0}^{T}t^i\|\nabla\dot{ u }\|_{L^2}^{2}dt\leq C,
\end{equation}
and
\begin{equation}\label{i313}
\sup_{t\in[0,T]}t^i\left(\|\na^2 u\|_{L^2}^2+\|\na P\|_{L^2}^2\right) \leq C.
\end{equation}
\end{lemma}

{\it Proof.}  First,
operating $\partial_{t}+u \cdot\nabla$ to \eqref{1.1}$_2^j$ yields that
\begin{eqnarray*}
\partial_{t}(\rho\dot{u}^j)+\div(\rho u \dot{u}^j)-\mu\Delta\dot{u}^j
= -\mu\partial_{i}(\partial_{i} u \cdot\nabla u^{j})
-\mu\div({\partial_{i} u \partial_{i}u^j})
-\partial_{j}\partial_{t}P-( u \cdot\nabla)\partial_{j}P.
\end{eqnarray*}

Now, multiplying the above equality by $\dot{u}^j$, we obtain after integration by parts that
\begin{align}\label{3.15}&
\frac12\frac{d}{dt}\int\rho|\dot{ u }|^{2}dx+\mu\int|\nabla\dot{ u }|^{2}dx \notag \\
&=-\int\mu\partial_{i}(\partial_{i} u \cdot\nabla u^{j})\dot{u}^jdx
-\int\mu\div({\partial_{i} u \partial_{i}u^j})\dot{u}^jdx+J \notag \\ &
\leq C \|\nabla u\|_{L^4}^4+\frac{\mu}{4}\|\nabla\dot{u}\|_{L^2}^2+J,
\end{align}
where
$$ J\triangleq -\int\dot{u}^j \partial_{t}\partial_{j}Pdx -\int\dot{u}^j u \cdot\nabla \partial_{j}P dx $$
satisfies
\begin{align}\label{i41}
J & =\int P_{t} \div\dot{u}dx+\int u \cdot\nabla P \div\dot{u}dx
+\int \partial_{j}u^{i}\partial_{i} P\dot{u}^{j}dx \notag \\
& = \int (P_{t}+ u \cdot\nabla P)\partial_{j}u^{i}\partial_{i}u^{j}dx-\int P\partial_{j}u^{i}\partial_{i}\dot{u}^{j}dx
\end{align}
due to $\div u=0.$

Then, since $\div u=0$, integration by parts implies
\begin{align*}
&\int (P_{t}+ u \cdot\nabla P)\partial_{j}u^{i}\partial_{i}u^{j}dx\\&=\frac{d}{dt}\int P\partial_{j}u^{i}\partial_{i}u^{j}dx-\int P\partial_{j}u_{t}^{i}\partial_{i}u^{j}dx-\int P\partial_{j}u^{i}\partial_{i}u_{t}^{j}dx\\
&\quad- \int P u^k\p_k(\p_ju^i\p_iu^j)dx\\
&=\frac{d}{dt}\int P\partial_{j}u^{i}\partial_{i}u^{j}dx -\int P\partial_{j}u_{t}^{i}\partial_{i}u^{j}dx-\int P u^k\p_k\p_ju^i\p_iu^jdx\\
&\quad -\int P\partial_{j}u^{i}\partial_{i}u_{t}^{j}dx-\int P u^k\p_k\p_iu^j\p_ju^idx\\
&=\frac{d}{dt}\int P\partial_{j}u^{i}\partial_{i}u^{j}dx  -\int P\partial_{j}\dot u ^{i}\partial_{i}u^{j}dx+\int P \p_j u^k\p_ku^i\p_iu^jdx\\
&\quad-\int P\partial_{j}u^{i}\partial_{i}\dot u^{j}dx+\int P \p_iu^k\p_ku^j\p_ju^idx,
\end{align*}
which together with \eqref{i41}, H{\"o}lder's and Young's inequalities yields
\begin{align}\label{i46}
J \leq & \frac{d}{dt}\int P\partial_{j}u^{i}\partial_{i}u^{j}dx
+C\int|P||\nabla\dot{ u }||\nabla u |dx
+C\int |P||\nabla u |^{3}dx \notag \\
\leq & \frac{d}{dt}\int P\partial_{j}u^{i}\partial_{i}u^{j}dx
+C \left(\|P\|_{L^4}^{4}+\|\nabla u \|_{L^4}^{4}\right)
+\frac{\mu}{4}\|\nabla\dot{ u }\|_{L^2}^{2}.
\end{align}

Next, putting \eqref{i46} into \eqref{3.15} gives
\be\ba\label{i47}
\Psi'(t)+\frac{\mu}{2}\int|\nabla\dot{ u }|^{2}dx\le C\|P\|_{L^4}^{4}+C\|\nabla u
\|_{L^4}^{4},
\ea\ee
where
\begin{equation*}
\Psi(t)\triangleq\frac12\int\rho|\dot{u}|^{2}dx
-\int P\partial_{j}u^{i}\partial_{i}u^{j}dx
\end{equation*}
satisfies \be\la{psi1} \frac14\int\n |\dot u|^2dx-C\|\na u\|_{L^2}^4\le \Psi(t)\le \int\n |\dot u|^2dx+C\|\na u\|_{L^2}^4\ee due to \eqref{3.8} and \eqref{001}.
Moreover, it follows from Sobolev's inequality, \eqref{001}, and \eqref{3.1} that
\begin{equation}\label{l1}\ba
\|P\|_{L^4}^{4}+\|\nabla u \|_{L^4}^{4}&\leq C\left(\|\nabla P\|_{L^{4/3}}^{4}+ \|\nabla^2u\|_{L^{4/3}}^{4}\right)\\&\leq C\|\rho\dot{ u }\|_{L^{4/3}}^{4}
\leq C\|\rho\|_{L^2}^{2}\|\sqrt{\rho}\dot{ u }\|_{L^2}^{4}
\leq C\|\sqrt{\rho}\dot{ u }\|_{L^2}^{4}.\ea
\end{equation}
Multiplying \eqref{i47} by $t^i (i=1,2)$ and using \eqref{psi1}--\eqref{l1},  we obtain \eqref{133} from Gronwall's inequality, \eqref{3.5}, and \eqref{022}.

Finally, \eqref{i313} is a direct consequence of \eqref{133} and \eqref{001}. The proof of Lemma \ref{lem3.4} is finished. \hfill $\Box$

\subsection{Higher Order Estimates}
The following spatial weighted estimate on the density plays an important role in deriving the bounds on the higher order derivatives of the solutions $(\n,u,P)$.

\begin{lemma}\label{lem03.6}
There exists a positive constant $C$ depending only on $\mu$, $\|\n_0\|_{L^1\cap L^\infty}$, $\|\n_0\bar x^a\|_{L^1}$, $\|\sqrt{\n_0} u_0\|_{L^2}$, $\|\na u_0\|_{L^2}$, $N_0$, and $T$ such that
\begin{equation}\label{06.1}
\sup_{t\in[0,T]}\|\rho\bar{x}^{a}\|_{L^{1}}\leq C.
\end{equation}
\end{lemma}

{\it Proof.}
First, for $N>1,$ let $\vp_N\in C^\infty_0(\r^2)$  satisfy
\be\ba \la{vp1}
0\le \vp_N \le 1, \quad  \vp_N(x)
=\begin{cases} 1,~~~~ |x|\le N/2,\\
0,~~~~ |x|\ge N,\end{cases}
\quad |\na \vp_N|\le C N^{-1}.
\ea\ee
It follows from \eqref{1.1}$_1$ that
\begin{align}\la{oo0}
\frac{d}{dt}\int \n \vp_{N} dx &=\int \n u \cdot\na \vp_{N} dx \notag \\
&\ge - C N^{-1}\left(\int\n dx\right)^{1/2}\left(\int\n |u|^2dx\right)^{1/2}\ge - \ti C N^{-1},
\end{align}
where in the last inequality one has used \eqref{3.1} and \eqref{0.11}.
Integrating \eqref{oo0} and choosing $N=N_1\triangleq2N_0+4\tilde CT$, we obtain after using \eqref{13} that
\begin{align}\la{p1}
\inf\limits_{0\le t\le T}\int_{B_{N_1}} \n dx&\ge \inf\limits_{0\le t\le T}\int \n \vp_{N_1} dx \notag \\
&\ge \int \n_0 \vp_{N_1} dx-\ti CN_1^{-1}T \notag \\
&\ge \int_{B_{N_0}} \n_0 dx-\frac{\ti C T}{2N_0+4\tilde C T} \notag \\
&\ge 1/4.
\end{align}
Hence, it follows from \eqref{p1}, \eqref{3.1}, \eqref{22}, \eqref{0.11}, and \eqref{3.5} that for any $\eta\in(0,1]$ and any $s>2$,
\begin{equation}\label{06.2}
\| u \bar{x}^{-\eta}\|_{L^{s/\eta}}\leq C\left(\|\n^{1/2}u\|_{L^2}+\|\na u\|_{L^2}\right)\le C.
\end{equation}
Multiplying \eqref{1.1}$_1$ by $\bar{x}^{a}$ and integrating the resulting equality by parts over $\mathbb{R}^2$ yield that
\bnn\ba
\frac{d}{dt}\int\rho\bar{x}^{a}dx & \leq C\int\rho| u |\bar{x}^{a-1}\log^{2}(e+|x|^2)dx\\
 &  \leq C\|\rho\bar{x}^{a-1+\frac{8}{8+a}}\|_{L^{\frac{8+a}{7+a}}}\| u \bar{x}^{-\frac{4}{8+a}}\|_{L^{8+a}} \\
 &  \leq C\int\rho\bar{x}^{a}dx+C,
\ea\enn
which along with Gronwall's inequality gives \eqref{06.1} and finishes the proof of Lemma \ref{lem03.6}.    \hfill $\Box$

\begin{lemma}\label{lem3.5}
There exists a positive constant $C$ depending on $T$ such that
\begin{align}\label{16.1}
\sup_{t\in[0,T]}\|\rho\|_{H^{1}\cap W^{1,q}}
&+\int_{0}^{T}\left(\|\nabla^{2} u \|_{L^2}^2+\|\nabla^{2} u \|_{L^q}^{\frac{q+1}{q}}
+t\|\nabla^{2} u \|_{L^2 \cap L^q}^2\right)dt \notag \\
&+\int_{0}^{T}\left(\|\nabla P\|_{L^2}^2+\|\nabla P \|_{L^q}^\frac{q+1}{q}+t\|\nabla P\|_{L^2\cap L^q}^2\right)dt
\leq C(T).
\end{align}
\end{lemma}

{\it Proof.}
First, it follows from the mass equation \eqref{1.1}$_1$ that $\na\n$ satisfies for any $r\ge 2,$
\begin{equation}\label{t}
\frac{d}{dt}\|\nabla\rho\|_{L^r}\leq C(r)\|\nabla u \|_{L^\infty}\|\nabla\rho\|_{L^r}.
\end{equation}
Next, one gets from Gagliardo-Nirenberg inequality, \eqref{3.5}, and \eqref{001} that for $q>2$ as in Theorem \ref{thm1},
\be\ba \label{4.5}
\|\nabla u \|_{L^\infty}  &  \leq C\|\nabla u\|_{L^2}^{\frac{q-2}{2(q-1)}}\|\nabla^2 u\|_{L^q}^{\frac{q}{2(q-1)}}  \leq C\|\rho\dot{ u }\|_{L^q}^{\frac{q}{2(q-1)}} .
\ea\ee

It follows from \eqref{p1}, \eqref{3.1}, \eqref{22}, and \eqref{06.1} that for any $\eta\in(0,1]$ and any $s>2$,
\begin{align}\label{lvb01}
\|\rho^\eta v\|_{L^{\frac{s}{\eta}}}
 &  \leq C\|\rho^\eta\bar x^{\frac{3\eta a}{4s}}\|_{L^{\frac{4s}{3\eta}}}
\|v\bar  x^{-\frac{3\eta a}{4s}}\|_{L^{\frac{4s}{\eta}}} \notag \\
 &  \leq C\|\n\|_{L^\infty}^{\frac{(4s-3)\eta}{4s}}\|\n\bar x^a\|_{L^1}^{\frac{3\eta}{4s}}\left( \|\n^{1/2} v\|_{L^2}+\|\na v\|_{L^2}\right) \notag \\
 &  \leq C\left(\|\n^{1/2}v\|_{L^2}+\|\na v\|_{L^2}\right),
\end{align}
which together with the Gagliardo-Nirenberg inequality shows that
\begin{align}\label{4.6}
\|\rho\dot{ u }\|_{L^q}
 &  \leq C\|\rho\dot{ u }\|_{L^2}^{\frac{2(q-1)}{q^2-2}}
\|\rho\dot{ u }\|_{L^{q^2}}^{\frac{q(q-2)}{q^2-2}} \notag \\
 &  \leq C\|\rho\dot{ u }\|_{L^2}^{\frac{2(q-1)}{q^2-2}}
\left(\|\sqrt{\rho}\dot{ u }\|_{L^2}+\|\nabla\dot{ u }\|_{L^2}\right)^{\frac{q(q-2)}{q^2-2}} \notag \\
 &  \leq C\left(\|\sqrt{\rho}\dot{ u }\|_{L^2}
+\|\sqrt{\rho}\dot{ u }\|_{L^2}^{\frac{2(q-1)}{q^2-2}}
\|\nabla\dot{ u }\|_{L^2}^{\frac{q(q-2)}{q^2-2}}\right).
\end{align}
This combined with \eqref{3.5} and \eqref{133} leads to
\begin{align}\label{4.7}
 &  \int_{0}^{T}\left(\|\rho\dot{ u }\|_{L^q}^{\frac{q+1}{q}}
    +t\|\rho\dot{ u }\|_{L^q}^{2}\right)dt \notag \\
 &  \leq C\int_{0}^{T}\left(\|\sqrt{\rho}\dot{ u }\|_{L^2}^{2}
    +t\|\nabla\dot{ u }\|_{L^2}^{2}
    +t^{-\frac{q^{3}-q^{2}-2q-1}{q^{3}-q^{2}-2q}}+1\right)dt \notag \\
 &  \leq C,
\end{align}
which along with \eqref{4.5} in particular implies
\be \la{nu1}\int_0^T\|\na u\|_{L^\infty}dt\le C. \ee
Thus, applying Gronwall's inequality to \eqref{t} gives
\begin{equation}\label{4.8}
\sup_{t\in[0,T]}\|\nabla\rho\|_{L^2\cap L^q}\leq C.
\end{equation}

Furthermore, it is easy to deduce from \eqref{001}, \eqref{4.7}, \eqref{3.1}, \eqref{3.5}, and \eqref{022} that
\begin{align}\label{4.10}
& \int_{0}^{T}\left(\|\nabla^2 u \|_{L^2}^2+\|\nabla^2 u \|_{L^q}^\frac{q+1}{q}+t\|\nabla^2 u \|_{ L^2 \cap L^q}^2\right)dt \notag \\
& +\int_0^T\left(\|\nabla P\|_{L^2}^2+\|\nabla P \|_{L^q}^\frac{q+1}{q}+t\|\nabla P\|_{L^2\cap L^q}^2\right)dt\leq C,
\end{align}
which together with \eqref{3.1} and \eqref{4.8} yields \eqref{16.1} and completes the proof of Lemma \ref{lem3.5}.  \hfill $\Box$

\begin{lemma}\label{lem3.6}
There exists a positive constant $C$ depending on $T$ such that
\begin{equation}\label{6.1}
\sup_{t\in[0,T]}\|\rho\bar{x}^{a}\|_{L^1\cap H^{1}\cap W^{1,q}}\leq C(T).
\end{equation}
\end{lemma}

{\it Proof.}
One derives from \eqref{1.1}$_1$ that $\rho\bar{x}^a$ satisfies
\begin{equation}\label{A}
\partial_{t}(\rho\bar{x}^a)+ u \cdot\nabla(\rho\bar{x}^a)
-a\rho\bar{x}^{a} u \cdot\nabla\log\bar{x}=0.
\end{equation}
Taking the $x_i$-derivative on the both side of \eqref{A} gives
\begin{align}\label{B}
0=  &  \partial_{t}\partial_{i}(\rho\bar{x}^a)+ u \cdot\nabla\partial_{i}(\rho\bar{x}^a)
+\partial_{i} u \cdot\nabla(\rho\bar{x}^a)
-a\partial_{i}(\rho\bar{x}^{a}) u \cdot\nabla\log\bar{x} \notag \\
 &  -a\rho\bar{x}^{a}\partial_{i} u \cdot\nabla\log\bar{x}
-a\rho\bar{x}^{a} u \cdot\partial_{i}\nabla\log\bar{x}.
\end{align}
For any $r\in[2,q]$, multiplying \eqref{B} by $|\na (\rho\bar{x}^a)|^{r-2}\partial_{i}(\rho\bar{x}^a)$ and integrating the resulting equality over $\mathbb{R}^2$, we obtain after integration by parts that
\begin{align}\label{6.4}
\frac{d}{dt}\|\nabla(\rho\bar{x}^a)\|_{L^r} \leq  &
C\left(1+\|\nabla u \|_{L^\infty}+\| u
\cdot\nabla\log\bar{x}\|_{L^\infty}\right)
\|\nabla(\rho\bar{x}^a)\|_{L^r} \notag \\
&  + C\|\rho\bar{x}^a\|_{L^\infty}\left(\||\nabla u ||\nabla\log\bar{x}|\|_{L^r}
+\||u||\nabla^{2}\log\bar{x}|\|_{L^r}\right) \notag \\
\leq & C\left(1+\|\nabla u\|_{W^{1,q}}\right)\|\nabla(\rho\bar{x}^a)\|_{L^r}  \notag \\
& +C\|\rho\bar{x}^a\|_{L^\infty}\left(\|\nabla u \|_{L^r}
+\|u \bar{x}^{-\frac{2}{5}}\|_{L^{4r}}\|\bar{x}^{-\frac{3}{2}}
\|_{L^{\frac{4r}{3}}}\right) \notag \\
\leq & C\left(1+\|\nabla u \|_{W^{1,q}}\right)
\left(1+\|\nabla(\rho\bar{x}^a)\|_{L^r}+\|\nabla(\rho\bar{x}^a)\|_{L^q}\right),
\end{align}
where in the second and the last inequalities, one has used \eqref{06.2} and \eqref{06.1}, respectively.
Choosing $r=q$ in \eqref{6.4}, together with \eqref{16.1}, indicates that
\begin{equation}\label{6.5}
\sup_{t\in[0,T]}\|\nabla(\rho\bar{x}^a)\|_{L^q}\leq C.
\end{equation}
Setting $r=2$ in \eqref{6.4}, we deduce from \eqref{16.1} and \eqref{6.5} that
\begin{equation}\label{14.2}
\sup_{t\in[0,T]}\|\nabla(\rho\bar{x}^a)\|_{L^2}\leq C.
\end{equation}
This combined with \eqref{6.5} and \eqref{06.1} gives \eqref{6.1} and finishes the proof of Lemma \ref{lem3.6}.   \hfill $\Box$

\begin{lemma}\label{lem3.8}
There exists a positive constant $C$ depending on $T$ such that
\begin{equation}\label{7.1}
\sup_{t\in[0,T]}t\|\sqrt{\rho}u_t \|_{L^2}^{2}
+\int_{0}^{T}t\|\nabla u_t \|_{L^2}^2dt\leq C(T).
\end{equation}
\end{lemma}

{\it Proof.}
Differentiating \eqref{1.1}$_2$ with respect to $t$ leads to
\begin{equation}\label{7.2}
\rho u_{tt}+\rho u \cdot\nabla u_t-\mu\Delta u_t
+\nabla P_t=-\rho_{t} u_t-(\rho u )_{t}\cdot\nabla u .
\end{equation}
Multiplying \eqref{7.2} by $u_t$ and integrating the resulting equation by parts over $\mathbb{R}^2$,
we obtain after using \eqref{1.1}$_1$ and \eqref{1.1}$_3$ that
\begin{align}\label{7.3}
\frac{1}{2}\frac{d}{dt}\|\sqrt{\rho} u_t \|_{L^2}^{2}
+\mu\|\nabla u_t\|_{L^2}^2
& =-\int \rho_{t}|u_t |^2dx
-\int (\rho u )_{t}\cdot\nabla u \cdot u_t dx \notag \\
&  \leq C\int\rho| u ||u_t |(|\nabla u_t|+|\nabla u |^2
+| u ||\nabla^2 u |)dx \notag \\
& \quad +C\int\rho| u |^2|\nabla u||\nabla u_t|dx+C\int\rho|u_t|^2|\nabla u |dx \notag \\
& \triangleq K_1+K_2+K_3.
\end{align}

We estimate each term on the right-hand side of \eqref{7.3} as follows.

First, the combination of \eqref{06.2}, \eqref{lvb01}, \eqref{0.11}, and \eqref{3.5} gives that for any
$\eta \in(0,1]$ and any $s>2$,
\begin{equation}\label{11.1}
\|\rho^{\eta} u\|_{L^{s/\eta}}+\|u \bar{x}^{-\eta}\|_{L^{s/\eta}}\leq C,
\end{equation}
which together with \eqref{lvb01}, \eqref{3.5}, and H{\"o}lder's inequality yields that
\begin{align}\label{7.5}
 K_1 &  \leq C\|\sqrt{\rho} u \|_{L^6}\|\sqrt{\rho}u_t\|_{L^2}^{1/2}
\|\sqrt{\rho}u_t\|_{L^6}^{1/2}
\left(\|\nabla u_t\|_{L^2}+\|\nabla u \|_{L^4}^2\right) \notag \\
 & \quad +C\|\rho^{1/4} u \|_{L^{12}}^{2}\|\sqrt{\rho}u_t\|_{L^2}^{1/2}
\|\sqrt{\rho}u_t \|_{L^6}^{1/2}\|\nabla^2 u \|_{L^2} \notag \\
 &  \leq C\|\sqrt{\rho}u_t\|_{L^2}^{1/2}
\left(\|\sqrt{\rho}u_t\|_{L^2}+\|\nabla u_t\|_{L^2}\right)^{1/2}
\left(\|\nabla u_t\|_{L^2}+\|\nabla^2 u \|_{L^2}\right)
\notag \\
&  \leq  \frac{\mu}{4}\|\nabla u_t\|_{L^2}^{2}+C
\left(1+\|\sqrt{\rho}u_t\|_{L^2}^{2}
+\|\nabla^2 u \|_{L^2}^2\right).
\end{align}
Next, H{\"o}lder's inequality combined with \eqref{lvb01}, \eqref{3.5}, and \eqref{11.1} leads to
\begin{align}\label{7.6}
K_2+K_3
&  \leq C\|\sqrt{\rho} u \|_{L^8}^{2}\|\nabla u \|_{L^4}\|\nabla u_t\|_{L^2}
+\|\nabla u \|_{L^2}\|\sqrt{\rho} u_t\|_{L^6}^{3/2}
\|\sqrt{\rho}u_t\|_{L^2}^{1/2} \notag \\
&  \leq \frac{\mu}{4}\|\nabla u_t\|_{L^2}^{2}
+C \left(1+\|\sqrt{\rho} u_t\|_{L^2}^{2}+\|\nabla^2 u \|_{L^2}^2\right).
\end{align}
Substituting \eqref{7.5}--\eqref{7.6} into \eqref{7.3} gives
\begin{equation}\label{11.3}
\frac{d}{dt}\|\sqrt{\rho}u_t\|_{L^2}^{2}
+\mu\|\nabla u_t\|_{L^2}^2\leq C\|\sqrt{\rho}u_t\|_{L^2}^{2}
+C\|\nabla^2 u \|_{L^2}^2+C.
\end{equation}

Finally, noticing that it follows from \eqref{3.5}, \eqref{11.1}, \eqref{16.1}, and Garliardo-Nirenberg inequality that
\begin{align}\label{jianew}
\int_0^T \|\sqrt{\rho} u_t\|_{L^2}^2dt&\le \int_0^T \left(\|\sqrt{\rho}\dot u\|_{L^2}^2+\|\sqrt{\rho}|u||\na u|\|_{L^2}^2\right)dt \notag \\
&\le C+C\int_0^T \|\sqrt{\rho} u \|_{L^6}^2\|\na u\|_{L^3}^2 dt \notag \\
&\le C+C\int_0^T \left(\|\na u\|_{L^2}^{2}+ \|\na^2 u\|_{L^2}^2\right) dt \notag \\
&\le C,
\end{align}
multiplying \eqref{11.3} by $t$, and using Gronwall's inequality, \eqref{16.1}, and \eqref{jianew}, we obtain \eqref{7.1} and complete the proof of Lemma \ref{lem3.8}. \hfill $\Box$

\section{Proof of Theorem \ref{thm1}}\label{sec4}
With the a priori estimates in Section \ref{sec3} at hand, we are now in a position to prove Theorem \ref{thm1}.

By Lemma 2.1, there exists a $T_{*}>0$ such that the problem \eqref{1.1}-\eqref{1.3} has a unique local strong solution $(\rho,u,P)$ on $\mathbb{R}^2\times(0,T_{*}]$. We plan to extend the local solution to all time.

Set
\begin{equation}\label{20.1}
T^{*}=\sup \{T~|~(\rho,u,P)\ \text{is a strong solution on}\ \mathbb{R}^{2}\times(0,T]\}.
\end{equation}
First, for any $0<\tau<T_*<T\leq T^{*}$ with $T$ finite, one deduces from \eqref{3.5}, \eqref{i313}, and
\eqref{7.1} that for any $q\geq2$,
\begin{equation}\label{20.2}
\nabla u \in C([\tau,T];L^2\cap L^q),
\end{equation}
where one has used the standard embedding
\begin{equation*}
L^{\infty}(\tau,T;H^1)\cap H^{1}(\tau,T;H^{-1})\hookrightarrow C(\tau,T;L^q)\ \ \text{for any}\ \ q\in[2,\infty).
\end{equation*}
Moreover, it follows from \eqref{16.1}, \eqref{6.1}, and \cite[Lemma 2.3]{L1996} that
\begin{equation}\label{20.3}
\rho\in C([0,T];L^1\cap H^1\cap W^{1,q}).
\end{equation}

Finally, if $T^{*}<\infty,$ it follows from \eqref{20.2}, \eqref{20.3}, \eqref{3.5}, \eqref{0.11}, and \eqref{6.1} that
$$(\n, u)(x,T^*)=\lim_{t\rightarrow T^*}(\n, u)(x,t)$$
satisfies the initial conditions \eqref{2.2} at $t=T^*$. Thus, taking $(\n, u)(x,T^*)$ as the initial data, Lemma 2.1 implies that one can extend the  strong solutions beyond $T^*$. This contradicts the assumption of $T^*$ in \eqref{20.1}. The proof of Theorem \ref{thm1} is completed.  \hfill $\Box$

\end{document}